\documentclass{amsart}
\usepackage{bbm}
\usepackage{hyperref}
\usepackage{dsfont}
\usepackage{mathrsfs,lscape}
\usepackage[all]{xy}
\input macros.tex

\begin{document}
\title{Advances in Losing}
\author{Thane E. Plambeck}
\maketitle
\begin{abstract}
We survey recent developments in the theory of impartial combinatorial games in misere play,
focusing on how the Sprague-Grundy theory
of normal-play impartial games generalizes to misere play via the {\em indistinguishability quotient
construction} \cite{ttw}.
This paper is based on a lecture given on 21 June 2005 at the Combinatorial Game Theory Workshop at the
Banff International Research Station.  It has been extended to include a survey of 
results on misere games, a list of open problems involving them, and a summary of {\em MisereSolver} \cite{miseresolver}, the excellent Java-language
program for misere indistinguishability quotient construction recently developed by Aaron Siegel.
Many wild misere games that have long appeared intractible may 
now lie within the grasp of assiduous losers and their
faithful computer assistants, particularly those researchers and computers equipped with {\em MisereSolver}.
\end{abstract}

\section{Introduction}
\begin{quotation}
{\sf We've spent a lot of time teaching you how to win games by being the last
to move.  But suppose you are baby-sitting little Jimmy and want, 
at least occasionally, to make sure you {\em lose}?  This means that
instead of playing the normal play rule in which whoever can't move is
the {\em loser}, you've switched to {\bf misere play} rule when
he's the {\em winner}.  Will this make much difference?  Not {\em always}...}
\end{quotation}

That's the first paragraph from the thirteenth
chapter (``Survival in the Lost World'') 
of Berlekamp, Conway, and Guy's 
encyclopedic work on combinatorial game theory, 
{\em Winning Ways for your Mathematical Plays} \cite{ww}.

And why ``not {\em always}?'' The misere analysis of an impartial combinatorial
game often proves to be far more difficult than it is in normal play.
To take a typical example, the normal play analysis
of {\bf Dawson's Chess} \cite{dawson} was
published as early as 1956 by Guy and Smith \cite{gs}, but even today, a complete
misere analysis hasn't been found\footnote{See \S\ref{dawson}.}.  Guy tells the story  \cite{guyhalmos}:

\begin{quotation}
{\sf
[Dawson's chess] \ is played on a $3 \times n$ board with white pawns
on the first rank and black pawns on the third.  It was posed as a {\em losing}
game (last-player-losing, now called {\bf misere}) so that capturing was obligatory.
Fortunately, (because we {\em still} don't know how to play misere Dawson's Chess)
I assumed, as a number of writers of that time and since have done, that the 
misere analysis required only a trivial adjustment of the normal (last-player-winning)
analysis.  This arises because Bouton, in his original analysis of Nim \cite{bouton},
had observed that only such a trivial adjustment was necessary to cover both normal and 
misere play...

But even for {\em impartial} games, in which the same options are available to both
players, regardless of whose turn it is to move, Grundy \& Smith \cite{grundysmith} showed that the general
situation in misere play soon gets very complicated, and Conway \cite{onag}, (p. 140)
confirmed that the situation can only be simplified to the microscopically small extent
noticed by Grundy \& Smith.

At first sight Dawson's Chess doesn't look like an impartial game, but if you know how pawns
move at Chess, it's easy to verify that it's equivalent to the game played with rows of
skittles in which, when it's your turn, you knock down any skittle, together with its
immediate neighbors, if any.
}
\end{quotation}

So misere play can be difficult.  But is it a hopeless situation?  It has often seemed so.
Returning to
chapter 13 in \cite{ww}, one encounters the {\em genus theory of impartial
misere disjunctive sums},
extended significantly from its original presentation in chapter 7 
(``How to Lose When You Must'') of Conway's
{\em On Numbers and Games}~\cite{onag}.  
But excluding the {\em tame games} that play like Nim in misere play, 
there's a remarkable paucity of example games that the 
genus theory completely resolves.  For example, the section 
``Misere Kayles'' from the 1982 first edition of \cite{ww} promises

\begin{quotation}
{\sf Although several tame games arise in Kayles (see Chapter 4), wild game's abounding
and we'll need all our [genus-theoretic] resources to tackle it...}
\end{quotation}

However, it turns out Kayles isn't ``tackled'' at all---after an extensive table of genus values 
to heap size 20, one finds the question

\begin{quotation}
{\sf Is there a larger single-row P-position?}
\end{quotation}

It was left to the amateur William L. Sibert \cite{sibert} to settle misere Kayles using
completely different methods.  One finds a description of his solution
at end of the updated Chapter 13 in the second edition of \cite{ww}, and also in \cite{sibert}.
In 2003, \cite{ww} summarized the situation as follows (pg 451):

\begin{quotation}
{\sf
Sibert's remarkable {\em tour de force} raises once again the question: are misere analyses really
so difficult?  A referee of a draft of the Sibert-Conway paper wrote ``the actual solution will have
no bearing on other problems,'' while another wrote ``the ideas are likely to be applicable to some
other games...''
}
\end{quotation}

\subsection{Misere play--the natural impartial game convention?}

When nonmathematicians play impartial games,
they tend to choose the misere play convention\footnote{``Indeed, if anything, misere
Nim is more commonly played than normal Nim...'' \cite{onag}, pg 136.}.  This was already recognized
by Bouton in his classic paper
``Nim, A Game with a Complete Mathematical Theory,''
~\cite{bouton}:

\begin{quotation}
{\sf
The game may be modified by agreeing that the player
who takes the last counter from the table {\em loses}.
This modification of the three pile [Nim] game seems to be more
widely known than that first described, but its theory is
not quite so simple...
}
\end{quotation}

But why do people prefer the misere play convention?  The answer may lie in Fraenkel's
observation that impartial games lack {\em boardfeel}, and simple
{\em Schadenfreude}\footnote{The joy we take in another's misfortune.}:

\begin{quotation}
{\sf For many MathGames, such as Nim, a player without prior knowledge of the strategy has no
inkling whether any given position is ``strong'' or ``weak'' for a player.  Even two positions
before ultimate defeat, the player sustaining it may be in the dark about the outcome, which 
will stump him.  The player has no boardfeel...} (\cite{frbib} pg. 3).
\end{quotation}

If both players are ``in the dark,'' perhaps it's only natural that the last player compelled to make
a move
in such a pointless game should be deemed the {\em loser}.  
Only when a mathematician gets involved are
things ever-so-subtly shifted toward the normal play convention, instead---but 
this is only because there is a 
simple and beautiful theory of normal-play impartial games---the Sprague-Grundy theory.  Secretly
computing nim-values, mathematicians win normal-play
impartial games time and time again.   Papers on normal play impartial games outnumber
misere play ones by a factor of perhaps fifty, or even more\footnote{Based on an informal count of
papers in the \cite{frbib} CGT bibliography.}.  

In the last twelve months it has become clear how to generalize such Sprague-Grundy nim-value computations
to misere play via {\em indistinguishability quotient construction} \cite{ttw}.  
As a result, many misere game problems that have long appeared intractible, or have
been passed over in silence as too difficult, have now been solved.  Still others, such as a Dawson's Chess, appear to remain
out of reach and await new ideas.  The remainder of this paper surveys this largely unexplored
territory.

\section{Two wild games}

We begin with two impartial games:  {\bf Pascal's Beans}---introduced here for the first time---and {\bf Guiles} (the octal game {\bf 0.15}).  Each has 
a relatively simple normal-play solution, but is {\em wild}\footnote{See Chapter 13 (``Survival in the Lost World'') in \cite{ww}, and the present paper's
section \ref{whatiswild} for more information on wild misere games.} 
in misere play.   Wild games are characterized by having misere play
that differs in an essential way\footnote{To be made precise in section \ref{whatiswild}.} from the play of misere Nim.  They often prove notoriously difficult to analyze completely.
Nevertheless, we'll give complete misere analyses for both Pascal's Beans and Guiles by using the key idea of the {\em misere indistinguishability quotient}, which
was first introduced in \cite{ttw}, and which we take up in earnest in section \ref{quotient}.

\section{Pascal's Beans}

\label{iqc}
{\bf Pascal's Beans} is a two-player impartial combinatorial game.  It's played with heaps of beans placed on Pascal's triangle, which is 
depicted in Figure
\ref{pascalstriangle}.
A legal move in the game is to slide a single bean either up
a single row and to the left one position, or alternatively up a single row and to the right one position in the triangle.  For example,
in Figure \ref{pascalstriangle}, a bean resting on the cell marked 20 could be moved to either cell labelled 10.

\begin{figure}[h]
\begin{eqnarray*}
\begin{array}{cccccccccccccccccc}
& & & & & & & & & 1 & & & & & & & &\\
& & & & & & & & 1 & & 1 & & & & & & &\\
& & & & & & & 1 & & 2 & & 1 & & & & & &\\
& & & & & & 1 & & 3 & & 3 & & 1 & & & & &\\
& & & & & 1 & & 4 & & 6 & & 4 & & 1 & & & &\\
& & & & 1 & & 5 & &10 & &10 & & 5 & & 1 & & &\\
& & & 1 & & 6 & &15 & & 20& &15 & & 6 & & 1 & &\\
& & 1 & & 7 & &21 & &35 & &35 & &21 & & 7 & & 1 &\\
& & & & &\vdots & & & & \vdots & & & & \vdots& & & & \\
\end{array}
\end{eqnarray*}
\caption{\label{pascalstriangle} The Pascal's Beans board.}
\end{figure}

The actual numbers in Pascal's triangle
are not relevant in the play of the game, except for the 1's that mark the non-interior, or 
``boundary'' positions of the board.  
In play of Pascal's Beans, a bean is considered out of play 
when it first reaches a boundary position of the triangle.
  The game ends when all beans have reached the boundary.

\subsection{Normal play}
In {\em normal play} of Pascal's Beans,
the last player to make a legal move is declared the {\em winner} of the game.  
Figure \ref{triangle-normal} shows the pattern of {\em nim values} that arises in the analysis of the game.
Using the figure, it's possible to quickly determine the best-play outcome of an arbitrary starting
position in Pascal's Beans using the {\em Sprague-Grundy theory} \ and the
{\em nim addition} operation $\oplus$.  
Provided one knows the ${\mathbb Z}_2 \times {\mathbb Z}_2$ addition table in Figure \ref{z2z2}, all is well---the 
{\em P-positions} (second-player winning positions) are precisely those that have nim value zero (ie, $\ast 0$), and every other position is an {\em N-position}
(or next-player win), of nim value $\ast 1$, $\ast 2$, or $\ast 3$.

{\small
\begin{figure}[h]

\begin{eqnarray*}
\begin{array}{cccccccccccccccccc}
& & & & & & & & & \ast 0 & & & & & & & &\\
& & & & & & & & \ast 0 & & \ast 0 & & & & & & &\\
& & & & & & & \ast 0 & & \ast 1 & & \ast 0 & & & & & &\\
& & & & & & \ast 0 & & \ast 2 & & \ast 2 & & \ast 0 & & & & &\\
& & & & & \ast 0 & & \ast 1 & & \ast 0 & & \ast 1 & & \ast 0 & & & &\\
& & & & \ast 0 & & \ast 2 & & \ast 2 & & \ast 2 & & \ast 2 & & \ast 0 & & &\\
& & & \ast 0 & & \ast 1 & & \ast 0 & &  \underline{{\bf \ast 0}} & & \ast 0 & & \ast 1 & & \ast 0 & &\\
& & \ast 0 & & \ast 2 & & \ast 2 & & \underline{{\bf\ast 1}} & & \underline{{\bf \ast 1}}& & \ast 2 & &  \ast 2 & & \ast 0 &\\
   & \ast 0 & & \ast 1 & & \ast 0 & & \underline{{\bf \ast 0}} & & \underline{{\bf \ast 0}}& & \underline{{\bf \ast 0}} & &  \ast 0 & & \ast 1 & & \ast 0\\
& & & & &\vdots & & & & \vdots & & & & \vdots& & & & \\
\end{array}
\end{eqnarray*}
\caption{\label{triangle-normal} The pattern of single-bean nim-values in normal play of Pascal's Beans.  Each interior value
is the {\em minimal excludant} (or {\em mex}) of the two nim values immediately above it.  The underlined entries
form the first three rows of an infinite subtriangle whose rows alternate between $\ast 0$ and $\ast 1$.}
\end{figure}
}

{\small
\begin{figure}[h]

\begin{eqnarray*}
\begin{array}{c|cccc}
\oplus & \ast 0 & \ast 1 & \ast 2 & \ast 3 \\ \hline
\ast 0 & \ast 0 & \ast 1 & \ast 2 & \ast 3 \\ 
\ast 1 & \ast 1 & \ast 0 & \ast 3 & \ast 2 \\ 
\ast 2 & \ast 2 & \ast 3 & \ast 0 & \ast 1 \\ 
\ast 3 & \ast 3 & \ast 2 & \ast 1 & \ast 0 
\end{array}
\end{eqnarray*}
\caption{\label{z2z2} Addition for normal play of Pascal's Beans.}
\end{figure}
}

\subsection{Misere play}
In {\em misere play} of Pascal's Beans, the last player to make a move is declared the {\em loser} of the game. 
Is it possible to give an analysis of misere Pascal's Beans that resembles the normal play analysis?  
The answer is yes---but the positions of the triangle can no longer be identified with nim heaps $\ast k$, 
and the rule for the misere addition is no longer given by nim addition.  Instead, both the values to
be identified with particular positions of the triangle and the desired misere addition
are given by a particular twelve-element commutative monoid $\mathcal M$, 
the {\em misere indistinguishability quotient}\footnote{See section \ref{quotient}.} of Pascal's Beans.  The monoid $\mathcal M$ 
has an identity 1 and is presentable 
using three generators and relations:

\[{\mathcal M} = \langle \ a, b, c \ | \ a^2=1,\ c^2 =1,\ b^3 = b^2c \ \rangle. \]

Assiduous readers might enjoy verifying that the identity $b^4=b^2$ follows from these relations,
and that a general word of the form $a^ib^jc^k$ ($i, j, k \geq 0$) will always reduce to one of the
twelve {\em canonical words}
\[ {\ \mathcal M} = \{1,\ a,\ b,\ ab,\ b^2,\ ab^2,\ c,\ ac,\ bc,\ b^2c,\ abc,\ ab^2c \}. \]

Amongst the twelve canonical words, three represent P-position types
\[ \mathcal P = \{ a,\ b^2, ac \}, \]
and the remaining nine represent N-position types:
\[ \mathcal N = \{1,\ b,\ ab,\ ab^2,\ c,\ bc,\ b^2c,\ abc,\ ab^2c \}. \]

Figure \ref{triangle-misere} shows the identification of positions of the triangle with elements
of $\mathcal M$.
{\small
\begin{figure}[h]
\begin{eqnarray*}
\begin{array}{cccccccccccccccccc}
& & & & & & & & & 1 & & & & & & & &\\
& & & & & & & & 1 & & 1 & & & & & & &\\
& & & & & & & 1 & & \framebox{a} & & 1 & & & & & &\\
& & & & & & 1 & & b & & b & & 1  & & & & & \\
& & & & & 1 & & a & & \framebox{$b^2$} & & a & & 1 & & & & \\
& & & & 1 & & b & & c & & c & & b & & 1 & & &\\
& & & 1 & & a & & b^2 & &  \underline{{\bf b^2}} & & b^2 & & a & & 1 & &\\
& & 1 & & b & & c & & \underline{{\bf ab^2}} & & \underline{{\bf ab^2}}& & c & &  b & & 1 &\\
   & 1 & & a & & b^2 & & \underline{{\bf b^2}} & & \underline{{\bf b^2}}& & \underline{{\bf b^2}} & &  b^2 & & a & & 1 \\
& & & & &\vdots & & & & \vdots & & & & \vdots& & & & \\
\end{array}
\end{eqnarray*}
\caption{\label{triangle-misere} Identifications for single-bean positions in misere play of Pascal's Beans. 
The values are elements of the misere indistinguishability quotient $\mathcal M$ of Pascal's Beans. The underlined entries
form the first three rows of an infinite subtriangle whose rows alternate between the two values $b^2$ and $ab^2$.}
\end{figure}
}

Although we've used multiplicative notation to represent the addition operation in the monoid $\mathcal M$,
we use it to analyze general misere-play Pascal's Beans positions 
just as we used the nim values of Figure \ref{triangle-normal} and nim addition in normal play.  For example,
suppose a Pascal's Beans position involves just two beans---one placed along the central axis of the triangle at 
each of the two boxed positions in Figure \ref{triangle-misere}.  Combining the corresponding entries $a$ and $b^2$ as monoid elements,
we obtain the element $ab^2$, which we've already asserted is an N-position.  What is the winning misere-play move?  From the lower bean,
at the position marked $b^2$, the only available moves are both to a cell marked $b$.  This move is of the form
\[ ab^2 \rightarrow ab, \]
ie, the result is another misere N-position type (ie, $ab$).  So this option is not a winning misere move.  But the cell marked $a$ has an available
move is to the boundary.   The resulting winning move is of the form
\[ ab^2 \rightarrow b^2, \]
ie, the result is $b^2$, a P-position type.


\section{Guiles}

\label{guiles}
Guiles can be played with heaps of beans.  The possible moves are to remove a heap of 1
or 2 beans completely, or to take two beans from a sufficiently large heap and partition
what is left into two smaller, nonempty heaps.  This is the octal game {\bf 0.15}.

\subsection{Normal play}
The nim values of the octal game {\bf Guiles} 
fall into a period 10 pattern.  See Figure \ref{guiles-normal}.
\begin{figure}[h]
\begin{center}
\begin{tabular}{c|cccccccccc}
  & 1 & 2 & 3 & 4 & 5 & 6 & 7 & 8 & 9 & 10 \\ \hline 
 0+  & $1$  & $1$  & $0$  & $1$  & $1$  &  $2$  & $2$  & $1$  & $2$  & $2$  \\  
 10+  & $1$  & $1$  & $0$  & $1$  & $1$  & 
$2$  & $2$  & $1$  & $2$  & $2$  \\  
 20+  & $1$  & $1$  & $0$  & $1$  & $1$  & $2$  & $2$  & $1$  & $2$  & $2$  \\  
 30+  & $1$  & $1$  & $0$  & $1$  & $1$ &  $\cdots$
&  &  &  & 
\end{tabular}
\end{center}
\caption{\label{guiles-normal}Nim values for normal play {\bf 0.15}}
\end{figure}
\subsection{Misere play}
Using his recently-developed Java-language computer program {\em MisereSolver}, 
Aaron Siegel \cite{ps} found that the misere indistinguishability quotient $\mathcal Q$
of misere Guiles is a (commutative) monoid of order 42.  It has the presentation

\begin{eqnarray*} {\mathcal Q} = \langle \ a,\ b,\ c,\ d,\ e,\ f,\ g,\ h,\ i 
& | &  a^2=1,\ b^4=b^2,\ bc=ab^3,\ c^2=b^2,\ b^2d=d,\ \\
&   &  cd=ad,\ d^3=ad^2,\ b^2e=b^3,\ de=bd,\ be^2=ace,\\ 
&   & ce^2=abe,\ e^4=e^2,\ bf=b^3,\ df=d,\ ef=ace,\ \\
&   & cf^2=cf,\ f^3=f^2,\ b^2g=b^3,\ cg=ab^3,\ dg=bd,\ \\
&   & eg=be,\ fg=b^3,\ g^2=bg,\ bh=bg,\ ch=ab^3,\ \\
&   & dh=bd,\ eh=bg,\ fh=b^3,\ gh=bg,\ h^2=b^2,\ \\ 
&   & bi=bg,\ ci=ab^3,\ di=bd,\ ei=be,\ fi=b^3,\ \\
&   & gi=bg,\ hi=b^2,\ i^2=b^2 \ \rangle.
\end{eqnarray*}

In Figure \ref{guiles-misere} we show the single-heap misere equivalences for Guiles.
It is a remarkable fact that this sequence is also periodic of length ten---it's just that
the (aperiodic) {\em preperiod} is longer (length 66), and a person needs to know the monoid 
${\mathcal Q}$!  The P-positions of Guiles are the precisely those positions equivalent to one of the words
\[ P = \{\ a,\ b^2,\ bd,\ d^2,\ ae,\ ae^2,\ ae^3,\ af,\ af^2,\ ag,\ ah,\ ai\ \}.\]

\begin{figure}[h]
\begin{center}
$\begin{array}{c|cccccccccc}
    & 1 & 2 & 3 & 4 & 5 & 6 & 7 & 8 & 9 & 10 \\ \hline 
0+  & a & a & 1 & a & a & b & b & a & b & b \\
10+ & a & a & 1 & c & c & b & b & d & b & e \\
20+ & c & c & f & c & c & b & g & d & h & i \\
30+ & ab^2 & abg & f & abg & abe & b^3 & h & d & h & h \\
40+ & ab^2 & abe & f^2 & abg & abg & b^3 & h & d & h & h \\
50+ & ab^2 & abg & f^2 & abg & abg & b^3 & b^3 & d & b^3 & b^3 \\
60+ & ab^2 & abg & f^2 & abg & abg & b^3 & b^3 & d & b^3 & b^3 \\
70+ & ab^2 & ab^2 & f^2 & ab^2 & ab^2 & b^3 & b^3 & d & b^3 & b^3 \\
80+ & ab^2 & ab^2 & f^2 & ab^2 & ab^2 & b^3 & b^3 & d & b^3 & b^3 \\
90+ & ab^2 & ab^2 & f^2 & ab^2 & ab^2 & b^3 & b^3 & d & b^3 & b^3 \\
100+ & 
\end{array}$
\end{center}
\caption{\label{guiles-misere} Misere equivalences for Guiles. }
\end{figure}

Knowledge of the monoid presentation ${\mathcal Q}$, its partition into N- and P-position types, and the single-heap
equivalences in Figure \ref{guiles-misere} suffices to quickly determine
the outcome of an arbitrary misere Guiles position.  For example, suppose
a position contains four heaps of sizes 4, 58, 68, and 78.  Looking up
monoid values in Figure \ref{guiles-misere}, we obtain the product

\begin{eqnarray*}
a \cdot d \cdot d \cdot d & = & ad^3 \\
                                  & = & a \cdot ad^2 \ \ (\mbox{relation } d^3=ad^2) \\                                 
                                  & = & d^2 \ \ (\mbox{relation } a^2=1)  \\                                
\end{eqnarray*}
We conclude that 4+58+68+78 is a misere Guiles P-position.

\section{The indistinguishability quotient construction}
\label{quotient}
What do these two solutions have in common?  They were both obtained via a computer program
called {\em MisereSolver}, by Aaron Siegel.  Underpinning {\em MisereSolver} is the notion of
the {\em indistinguishability quotient construction}.  Here, we'll sketch the main ideas of the indistinguishability
quotient construction only.  They are developed in detail in \cite{ttw}.  

Suppose ${\mathcal A}$ is a set of (normal, or alternatively, misere) impartial game positions that is closed under the
operations of game addition and taking options (ie, making moves).   Unless we say otherwise, we'll always be taking
${\mathcal A}$ to be the set of all positions that arise in the play of a specific game $\Gamma$, which we fix in advance.
For example, one might take
\begin{eqnarray*} 
\Gamma & = & \mbox{Normal-play Nim}, \\
{\mathcal A} & = & \mbox{All positions that arise in normal-play Nim}, 
\end{eqnarray*} 
or 
\begin{eqnarray*} 
\Gamma & = & \mbox{Misere-play Guiles}, \\
{\mathcal A} & = & \mbox{All positions that arise in misere-play Guiles}.
\end{eqnarray*}

Two games $G, H \in {\mathcal A}$ are then said to be {\em indistinguishable}, and we write the relation $G \ \rho \ H$, 
if for every game $X \in {\mathcal A}$,
the sums $G+X$ and $H+X$ have the same outcome (ie, are both N-positions, or are both P-positions).  Note in particular that
if $G$ and $H$ are indistinguishable, then they have the same outcome (choose $X$ to be the {\em endgame}---ie, the terminal position, with no options).

The indistinguishability relation $\rho$ is easily seen to be an equivalence relation on ${\mathcal A}$,
but in fact more is true---it's a {\em congruence} on ${\mathcal A}$ \cite{ttw}.  This follows because indistinguishability is {\em compatible} with addition;
ie, for every set of three games $G, H, X \in {\mathcal A}$:

\begin{equation}
G \ \rho \ H \implies (G + X) \ \rho \ (H +X). 
\end{equation} 

Now let's make the definition 
\[\rho G = \{\ H \in {\mathcal A} \ | \ G \ \rho \ H \}.\]

We'll call $\rho G$  the {\em congruence class of ${\mathcal A}$ modulo $\rho$ containing $G$}.
Because $\rho$ is a congruence, there is a well-defined
addition operation 
\[\rho G + \rho H = \rho (G+H) \]
on the set ${\mathcal A} / \rho$ of all congruence classes $\rho G$ of ${\mathcal A}$ modulo $\rho$
\begin{equation}
\label{qe}
{\mathcal Q} = {\mathcal Q}(\Gamma)  = {\mathcal A} / \rho = \{ \ \rho G \ | \ G \in {\mathcal A}. \ \}
\end{equation}
The monoid ${\mathcal Q}$  is called the {\em indistinguishability quotient} of $\Gamma$.  
It captures the essential information of ``how to add'' in the play of game $\Gamma$,
and is the central figure of our drama.

The natural mapping $\Phi$ from ${\mathcal A}$ to ${\mathcal A} / \rho$
\[\Phi: G \mapsto \rho G\]
is called a {\em pretending function} (see \cite{ttw}).  Figures \ref{triangle-misere} and \ref{guiles-misere} illustrate the (as it happens, provably 
periodic \cite{ttw}) pretending 
functions of Pascal's Beans and Guiles, respectively.  We shall gradually come to see that the recovery
of ${\mathcal Q}$ and $\Phi$ from $\Gamma$ is the essence of impartial combinatorial game analysis in both normal
and misere play.

When $\Gamma$ is chosen as a normal-play impartial game, the elements of ${\mathcal Q}$
work out to be in 1-1 correspondence with the {\em nim-heap
values} (or $G$-{\em values}) that occur in the play of the game $\Gamma$.   For if $G$ and $H$ are normal-play impartial games with
$G=\ast g$ and $H=\ast h$, one easily shows that $G$ and $H$ are indistinguishable if and only if $g=h$.
Additionally, in normal play, every position $G$ satisfies the equation
\[G+G=0.\] 
As a result, the addition in a normal-play indistinguishability quotient is an abelian group in which every element is its own additive inverse.
The addition operation in the quotient ${\mathcal Q}$ is {\em nim addition}.
Every normal play indistinguishability quotient is therefore isomorphic to a (possibly infinite) direct
product
\[{\mathbb Z}_2 \times {\mathbb Z}_2 \times \cdots, \]
and a position is a P-position precisely if it belongs the congruence class of the identity (ie, $\ast 0$) of this group.
In this sense ``nothing new'' is learned about normal play impartial games via the indistinguishability quotient 
construction---instead, we've simply recast the Sprague-Grundy theory in new language.  
The fun begins when the construction is applied
in {\em misere play}, instead.

\section{Misere indistinguishability quotients}

In misere play, the indistinguishability quotient ${\mathcal Q}$ turns out to be a commutative monoid whose structure intimately depends upon
the particular game $\Gamma$ that is chosen for analysis.  We need to cover some background
material first.

\subsection{Preliminaries}
Consider the following three concepts in impartial games:

\begin{enumerate}
\item The notion of the {\em endgame} (or {\em terminal position}), ie, a game that has no options at all.
\item The notion of a {\em P-position}, ie, a game that is a second-player win in best play of the game.
\item The notion of the {\em sum of two identical games}, ie $G+G$.
\end{enumerate}

In normal play, these three notions are {\em indistinguishable}---wherever a person sees (1) in a sum $S$, he could freely substitute
(2) or (3) (or vice-versa, or any combination of such substitutions) without changing the outcome of $S$.

The three notions do not concide in misere play.   Let's see what happens instead.

\subsubsection{The misere endgame} In misere play, the endgame
is an N-position, not a P-position:  even though there is no move available from the endgame, a player still wants it to be
his {\em turn to move} when facing the endgame in misere play, because that means his opponent just {\em lost}, on his previous move.

\subsubsection{Misere outcome calculation}  After the special case of the endgame is taken care of, the recursive rule for outcome calculation in misere
play is exactly as it is in normal play:  a non-endgame position $G$ is a P-position iff all its options are N-positions.    Misere games
cannot be identified with nim heaps, in general, however---instead, a typical misere game looks like a complicated, usually unsimplifiable tree of options
\cite{onag}, \cite{grundysmith}.

\subsubsection{Misere P-positions}
Since the endgame is not a misere P-position, the simplest misere P-position is
the {\em nim-heap of size one}, ie, the game played using one bean on a table, where the game is to take that bean.
To avoid confusion both with what happens in normal play, and with the algebra of the misere indistinguishability quotient to
be introduced in the sequel, let's introduce some special symbols for the three simplest misere games:
\begin{eqnarray*}
\mathbbm{o} & = & \mbox { The misere {\em endgame}, ie, a position with no moves at all.} \\
\mathbbm{1} & = & \mbox { The misere {\em nim heap of size one}, ie, a position with one move (to } \mathbbm{o}). \\
\mathbbm{2} & = & \mbox { The misere {\em nim heap of size two}, ie, the game } \{\mathbbm{o}, \mathbbm{1}\}.
\end{eqnarray*}

Two games that we've intentionally left off this list are $\{\mathbbm{1}\}$ and  $\mathbbm{1}+\mathbbm{1}$. 
Assiduous readers should verify they are both indistinguishable from $\mathbbm{o}$.

\subsubsection{Misere sums involving P-positions}
Suppose that $G$ is an arbitrary misere P-position. Consider the misere sum
\begin{equation}
\label{1plus}
S = \mathbbm{1}+ G.
\end{equation}

Who wins $S$?  It's an N-position---a winning first-player move is to simply take the nim heap of size one, leaving the opponent to move
first in the P-position $G$.  In terms of outcomes, equation (\ref{1plus}) looks like
\begin{equation}
\label{1out}
N = P + P.
\end{equation}
Equation (\ref{1out}) does not remind us of normal play very much---instead, we always have $P+P=P$ in normal play.
On the other hand, it's not true that sum of two misere P-positions is {\em always} a misere N-position---in fact, when two typical misere P-positions
$G$ and $H$ are added together with {\em neither} equal to $\mathbbm{1}$, it {\em usually} happens that their sum is a P-position,
also.  But that's not {\em always} the case---it's also possible that two misere impartial P-positions, neither of which is 
$\mathbbm{1}$, can nevertheless result in an N-position when added together.   Without knowing the details of the 
misere P-position involved, little more can be said in general about the outcome when it's added to another game.

\subsubsection{Misere sums of the form $G+G$}
In normal play, a sum $G+G$ of two identical games is always indistinguishable from the endgame.  In misere play,
it's true that both $\mathbbm{o} +\mathbbm{o}$ and $\mathbbm{1} +\mathbbm{1}$ are indistinguishable from $\mathbbm{o}$, but beyond
those two sums, positions of the form $G+G$ are rarely indistinguishable from $\mathbbm{o}$.  It frequently happens
that a position $G$ in the play of a game $\Gamma$ has no $H \in {\mathcal A}$ such that $G+H$ is indistinguishable
from $\mathbbm{o}$.  This lack of natural inverse elements makes the structure of a typical misere indistinguishability quotient  
a {\em commutative monoid} rather than an {\em abelian group}. 

\subsubsection{The game $\mathbbm{2}+ \mathbbm{2}$}
The sum
\[\mathbbm{2} + \mathbbm{2} \]
is an important one in the theory of impartial misere games.
It's a P-position in misere play:  for if you move first by taking 1 bean from one summand,
I'll take two from the other, forcing you to take the last bean.  Similarly, if you choose to take 2 beans, I'll take 1 from the other.
So whereas in normal play one has the equation
\[ (\ast 2 + \ast 2) \ \ \rho \  \ast 0,\]
it's certainly {\bf not} the case in misere play that
\[ (\mathbbm{2} + \mathbbm{2}) \ \ \rho \ \ \mathbbm{o}, \]
since the two sides of that proposed indistinguishability relation don't even have the same outcome.
But perhaps
\begin{equation}
\label{2plus2vs1}
\mathbbm{2} + \mathbbm{2} \ \stackrel{\mbox ?}{\rho} \ \mathbbm{1} \ \  
\end{equation}
is valid?  The indistinguishability relation (\ref{2plus2vs1}) looks plausible at first glance---at least the positions on both sides are P-positions.
To decide whether it's possible to distinguish between $\mathbbm{2} + \mathbbm{2}$ and $\mathbbm{1}$, we might try adding various fixed
games $X$ to both, and see if we ever get differing outcomes:

\begin{center}
\begin{tabular}{ccc}
Misere & Misere & Misere \\ 
game & outcome of & outcome of \\
$X$ & $\mathbbm{2}+\mathbbm{2}+X$ & $\mathbbm{1}+X$ \\ \hline
$\mathbbm{o}$ & P & P \\
$\mathbbm{1}$ & N & N \\
$\mathbbm{2}$ & N & N \\
$\mathbbm{1}+\mathbbm{2}$ & N & N \\
$\mathbbm{2}+\mathbbm{2}$ & P & N \\
\end{tabular}
\end{center}

The two positions look like they {\em might} be indistinguishable, until we reach the final row of the table.  It
reveals that $(\mathbbm{2} + \mathbbm{2})$ distinguishes between $(\mathbbm{2} + \mathbbm{2})$  and $\mathbbm{1}$.
So equation (\ref{2plus2vs1}) fails.
Since a set of misere game positions ${\mathcal A}$ that includes $\mathbbm{2}$ and is closed under addition and taking options 
must contain all of the games $\mathbbm{1}$, $\mathbbm{2}$, and $\mathbbm{2}+\mathbbm{2}$, we've shown that a game that isn't
She-Loves-Me-She-Loves-Me-Not {\em always} has at least {\em two} distinguishable P-position types.  In normal play, there's just one P-position
type up to indistinguishability---the game $\ast 0$.

\subsection{Indistinguishability vs canonical forms}
In normal play, the Sprague-Grundy theory describes how to determine the outcome of a sum $G+H$ of two games $G$
and $H$ by computing {\em canonical} (or {\em simplest}) forms for each summand---these turn out to be {\em nim-heap 
equivalents} $\ast k$.   In both normal and misere play, canonical forms are obtained by pruning reversible moves from game trees
(see \cite{grundysmith}, \cite{onag} and \cite{ww}). 

In \cite{onag}, Conway succinctly gives the rules for misere game tree simplification to canonical form:

{\sf
\begin{quotation}
When $H$ occurs in some sum we should naturally like to replace it by [a] simpler game $G$.  Of course,
we will normally be given only $H$, and have to find the simpler game $G$ for ourselves.  How do we do this?
Here are two observations which make this fairly easy:
\begin{enumerate}
\item $G$ must be obtained by deleting certain options of $H$.
\item $G$ itself must be an option of any of the deleted options of $H$, and so $G$ must be itself be a {\em second option}
of $H$, if we can delete any option at all.
\end{enumerate}
On the other hand, if we obey (1) and (2), the deletion is permissible, except that we can only delete {\em all} the options
of $H$ (making $G$ = 0 [the endgame]) if one of the them is a second-player win.
\end{quotation}
}

Unlike in normal play, the canonical form of a misere game is not a nim heap in general.  
In fact, many misere game trees hardly simplify at all under the misere simplification rules.
Figure \ref{conway-list}, which duplicates information in \cite{onag} (its Figure 32), shows the 22 misere game trees born by day 4.

\begin{figure}[h]
\[\begin{array}{ccccc}
& \mbox{         } & & \mbox{        } & \\
\mathbbm{o}=\{\} & & \mathbbm{2}_{++} = \{\mathbbm{2_+}\}  & & \mathbbm{2}_{+}3\mathbbm{o} = \{\mathbbm{2}_+,3,\mathbbm{o}\} \\
\mathbbm{1}=\{\mathbbm{o}\} & & \mathbbm{2}_{+}\mathbbm{o} = \{\mathbbm{2}_+,\mathbbm{o}\} & & \mathbbm{2}_{+}3\mathbbm{1} = \{\mathbbm{2}_+,3,\mathbbm{1}\} \\
\mathbbm{2}=\{\mathbbm{o,1}\} & & \mathbbm{2}_{+}\mathbbm{1} = \{\mathbbm{2}_+,\mathbbm{1}\} & & \mathbbm{2}_{+}3\mathbbm{2} = \{\mathbbm{2}_+,3,\mathbbm{2}\} \\
3=\{\mathbbm{o,1,2}\} & & \mathbbm{2}_{+}\mathbbm{2} = \{\mathbbm{2}_+,\mathbbm{2}\} & & \mathbbm{2}_{+}3\mathbbm{2o} = \{\mathbbm{2}_+,3,\mathbbm{2,o}\} \\
4=\{\mathbbm{o,1,2},3\} & & \mathbbm{2}_{+}\mathbbm{2o} = \{\mathbbm{2}_+,\mathbbm{2,o}\} & & \mathbbm{2}_{+}3\mathbbm{21} = \{\mathbbm{2}_+,3,\mathbbm{2,1}\} \\
\mathbbm{2}_+ = \{\mathbbm{2}\} & & \mathbbm{2}_{+}\mathbbm{21} = \{\mathbbm{2}_+,\mathbbm{2,1}\} & & \mathbbm{2}_{+}3\mathbbm{21o} = \{\mathbbm{2}_+,3,\mathbbm{2,1,o}\} \\
3_+=\{3\} & & \mathbbm{2}_{+}\mathbbm{21o} = \{\mathbbm{2}_+,\mathbbm{2,1,o}\} & &  \\
\mathbbm{2}+\mathbbm{2}=\{3,\mathbbm{2}\} & & \mathbbm{2}_{+}3 = \{\mathbbm{2}_+,3\} & & \\
\end{array}\]
\caption{\label{conway-list} Canonical forms for misere games born by day 4.}
\end{figure}
Whereas only one normal-play nim-heap is born at each
birthday $n$, over 4 million nonisomorphic misere canonical forms are born by day five.  The number continues to grow very rapidly, roughly like a tower of exponentials of height $n$ (\cite{onag}).  
This very large
number of mutually distinguishable trees has often made misere analysis look like a hopeless activity.

\subsubsection{Indistinguishability identifies games with different misere canonical forms}
The key to the success of the indistinguishability quotient construction is that it is a {\em construction localized to the play of a particular 
game $\Gamma$}.  It
therefore has the possibility of identifying misere games
with different canonical forms.  While it's true that for misere games $G$, $H$ with different canonical forms
that there must be a game $X$ such that $G+X$ and $H+X$ have different outcomes, such an $X$ {\em might possibly
never occur} in play of the fixed game $\Gamma$ that we've chosen to analyze.  Indistinguishability
quotients are often {\em finite}, even for games $\Gamma$ that involve an infinity of different canonical forms 
amongst their position sums.

\section{What is a wild misere game?}
\label{whatiswild}
Roughly speaking, a misere impartial game $\Gamma$ is said to be {\em tame} when a complete analysis of it
can be given by identifying each of its positions with some position that arises in the misere play of Nim.  Tameness
is therefore an attribute of a {\em set} of positions, rather than a {\em particular} position.   Games $\Gamma$ that are
not tame are said to be {\em wild}.  Unlike tame games, wild games cannot be completely analyzed by viewing them as disguised versions of misere Nim.

\subsection{Tame games}

Conway's {\em genus theory} was first described in chapter 12 of \cite{onag}.
It describes a method for calculating whether all the positions of particular misere game $\Gamma$ are tame, and how to give a complete
analysis of $\Gamma$, if so.  For completeness, we've summarized the genus theory in the Appendix (section (\ref {gthm})) of this paper.

For misere games $\Gamma$
that the genus theory identifies as tame, a complete analysis can be given without reference to the 
indistinguishability quotient construction.  Various efforts to extend the genus theory to wider classes of games
have been made.   Example settings where progress has been made
are the main subject of papers by
of Ferguson \cite{ferguson2}, \cite{ferguson3} and Allemang \cite{a1}, \cite{a2}, \cite{a3}.  

\subsubsection{Indistinguishability quotients for tame games}

In this section, we reformulate the genus theory of tame games in terms of the indistinguishability quotient language.

Suppose $S$ is some finite set of misere combinatorial games.   We'll use the notation cl$(S)$ (the {\em closure} of $S$) to stand for
the smallest set of games that includes every element of $S$ and is closed under addition and taking options.  
Putting ${\mathcal A} = \mbox{cl}(S)$ and defining the indistinguishability quotient
\[{\mathcal Q} = A / \rho,\]
the natural question arises, what is the monoid ${\mathcal Q}$?  Figure \ref{simpleq} shows answers for $S = \{\mathbbm{1}\}$ and $S = \{\mathbbm{2}\}$.

\begin{figure}[h]
\[\begin{array}{ccccl}
  & \mbox{Presentation for} & & & \\
S & \mbox{monoid } {\mathcal Q} & \mbox{Order} & \mbox{Symbol} & \mbox{Name} \\ \hline
&  \\
\{\mathbbm{1}\} & \langle \ a \ | \ a^2 = 1 \ \rangle  & 2 & {\mathcal T}_1 &  \mbox{ First tame quotient} \\ \\
\{\mathbbm{2}\} & \langle \ a,\ b \ | \ a^2 = 1, \ b^3=b \ \rangle & 6 & {\mathcal T}_2 & \mbox{ Second tame quotient} \\

\end{array}\]
\caption{\label{simpleq} The first and second tame quotients }
\end{figure}

${\mathcal T}_1$ is called the {\em first tame quotient}.  It represents the misere play of {\em She-Loves-Me, She-Loves-Me-Not}.  
In ${\mathcal T}_1$, misere P-positions are represented by the 
monoid (in fact, group)
element $a$, and N-positions by 1.  

${\mathcal T}_2$, the {\em second tame quotient}, has the presentation
\[\langle\ a,b \ | \ a^2 = 1, \ b^3 = b \ \rangle. \]
It is a six-element monoid with two P-position types $\{a, b^2\}$.  The prototypical game $\Gamma$ with misere indistinguishability quotient $\mathcal T_2$ is the game of Nim, played with heaps of 1 and 2 only.
See Figures \ref{0333x} and \ref{corresp}.

\begin{center}
\begin{figure}[h]
\begin{picture}(45,120)(0,10)
\setlength{\unitlength}{.8\unitlength}


\put(5,10){$1$}
\put(14,20){\vector(1,1){30}}
\put(2,20){\vector(-1,1){30}}
\put(-24,54){\vector(1,-1){30}} 

\put(-35,55){$a$}
\put(-26,65){\vector(1,1){26}}

\put(45,55){$b$}
\put(44,66){\vector(-1,1){26}}
\put(55,64){\vector(1,1){26}}
\put(81,90){\vector(-1,-1){26}}

\put(2,95){$ab$}
\put(17,105){\vector(1,1){20}}
\put(37,125){\vector(-1,-1){20}}
\put(14,87){\vector(1,-1){26}} 

\put(85,95){$b^2$}
\put(85,107){\vector(-1,1){22}}
\put(59,125){\vector(1,-1){22}}

\put(42,130){$ab^2$}



\end{picture}
\caption{\label{0333x} {\em The misere impartial game theorist's coat of arms,} or the Cayley graph of $\mathcal T_2$.  Arrows have been drawn to show the 
action of the generators $a$ (the doubled rungs of the ladder) and $b$ 
(the southwest-to-northeast-oriented arrows) on $\mathcal{T}_{2}$.  See also Figure \ref{corresp}.}
\end{figure}
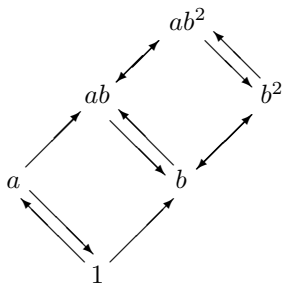
\end{center}

\begin{figure}[h]
\begin{tabular}{cccc}
 & Misere indistinguishability & &  \\ 
Position type & quotient element & Outcome & Genus \\ \hline \hline
& & & \\
Even \#$\mathbbm{1}$'s only &  1 &  N & $0^{120}$ \\ 
& & & \\ \hline
& & & \\ 
Odd \#$\mathbbm{1}$'s only &  $a$ &  P & $1^{031}$ \\  
& & & \\ \hline \\
Odd \#$\mathbbm{2}$'s &   &  &  \\
and & $b$ & N & $2^{20}$\\ 
Even \#$\mathbbm{1}$'s &   &  &  \\ \\ \hline \\
Odd \#$\mathbbm{2}$'s &   &  &  \\
and & $ab$ & N & $3^{31}$ \\
Odd \#$\mathbbm{1}$'s &   &  &  \\ \\ \hline \\

Even \#$\mathbbm{2}$'s ($\geq 2$) &   &  &  \\
and & $b^2$ & P & $0^{02}$ \\
Even \#$\mathbbm{1}$'s &   &  &  \\ \\ \hline \\

Even \#$\mathbbm{2}$'s ($\geq 2$) &   &  &  \\
and & $ab^2$ & N & $1^{13}$ \\
Odd \#$\mathbbm{1}$'s &   &  &  \\ \\ \hline \\

\end{tabular}
\caption{\label{corresp} When misere Nim is played with heaps of size 1 and 2 only, the resulting
misere indistinguishability quotient is the tame six-element monoid ${\mathcal T}_2$.  For more on genus
symbols and tameness, see section \ref{whatiswild}.  See also Figure \ref{0333x}.}
\end{figure}

\subsubsection{The general tame quotient}

For $n \geq 2$, the $n^{th}$ {\em tame quotient} is the monoid ${\mathcal T}_n$ with $2^n+2$ elements and the presentation

\begin{eqnarray*}
\begin{array}{lll}
{\mathcal T}_n = \langle\ a,\underbrace{b,c,d,e,f,g,\ldots}_{n-1 \mbox{ generators}} & | &  a^2=1, \\
                           &   & b^3=b, \ c^3=c, \ d^3=d, \ e^3=e, \ \ldots,  \\
                           &   & b^2=c^2=d^2=e^2= \ \ldots \ \rangle. 
\end{array} 
\end{eqnarray*}

${\mathcal T}_n$ is a disjoint union of its two maximal subgroups ${\mathcal T}_n = U \cup V.$
The set \[U = \{1, a\}\] is isomorphic to $\mathbb{Z}_2$.  The remaining $2^n$ elements of ${\mathcal T}_n$ form the set 
\begin{eqnarray*}
 V  =  \{ \ \ a^{a_i}b^{b_i}c^{c_i}d^{d_i}e^{e_i} \cdots  &  | & \ a_i = \mbox{0 or 1}   \\
&    & \ b_i = \mbox{1 or 2}  \\
&    & \ \mbox{Each of } \ c_i, d_i, e_i, \cdots = \mbox{0 or 1} \ \}.  \\
\end{eqnarray*}
and have an addition isomorphic to $\underbrace{\mathbb{Z}_2 \times \cdots \times \mathbb{Z}_2}_{n \mbox{ copies}}$.
The elements $a$ and $b^2$ are the only P-position types in ${\mathcal T}_n$.

\section{More wild quotients}

\subsection{The commutative monoid ${\mathcal R}_8$}  The smallest {\em wild} misere indistinguishability quotient 
${\mathcal R}_8$ has eight elements, and is unique up to isomorphism \cite{as1} amongst misere quotients with eight elements.
Its monoid presentation is

\[ {\mathcal R}_8 = \langle \ a,\ b,\ c\ \ | \ a^2=1,\ b^3=b,\ bc = ab,\ c^2=b^2 \ \rangle. \]
The P-positions are $\{a, b^2\}$.   

\subsubsection{{\bf 0.75}} 
An example game with misere
quotient ${\mathcal R}_8$ is the octal game {\bf 0.75}.  The first complete analysis of {\bf 0.75} was
given by Allemang 
using his {\em generalized genus theory} \cite{a1}.  Alternative formulations of the {\bf 0.75} solution are
also discussed at length in the appendix of \cite{plambeck} and in \cite{a2}.
See Figure \ref{075}.

\begin{figure}
\begin{center}
\[
\begin{array}{c|cc}
   & 1 & 2 \\ \hline
0+ & 1 & a \\
2+ & b & a \\
4+ & b & c \\
6+ & b & c \\
8+ & b & ab^2 \\
10+ & b & ab^2 \\
12+ & b & ab^2 \\
14+ & \ldots &  \\
\end{array}\]
\end{center}
\caption{\label{075} The pretending function for misere play of {\bf 0.75}.}
\end{figure}

\subsection{Flanigan's games} 

Jim Flanigan found solutions to the wild octal games {\bf 0.34} and {\bf 0.71};  a 
description of them
can be found in the ``Extras'' of chapter 13 in \cite{ww}.  It's interesting
to write down the corresponding misere quotients.  

\subsubsection{{\bf 0.34}}
The misere indistinguishability quotient of {\bf 0.34} has order 12.  There are three P-position types.  The pretending
function has period 8 (see Figure \ref{034}).

\[ {\mathcal Q}_{\bf 0.34} = \langle \ a, b, c \ | \ a^2=1,\ b^4=b^2,\ b^2c=b^3,\ c^2=1 \ \rangle \]

\[P = \{a,b^2,ac \}\]

\begin{figure}
\begin{center}
\[
\begin{array}{c|cccccccc}
   & 1 & 2 & 3 & 4 & 5 & 6 & 7 & 8 \\ \hline
0+ & a & 1 & a & b & 1 & a & 1 & ab \\
8+ & a & c & a & b & 1 & ac & 1 & ab \\
16+ & a & c & a & b & 1 & ac & 1 & ab 
\end{array}\]
\end{center}
\caption{\label{034} The pretending function for misere play of {\bf 0.34}.}
\end{figure}

\subsubsection{{\bf 0.71}}

The game {\bf 0.71} has a misere quotient of order 36 with the presentation

\[ \mathcal{Q}_{\bf 0.71} = \langle \ a,b,c,d \ | \ a^2=1, \ b^4=b^2, \ b^2c=c, \ c^4=ac^3, \ c^3d=c^3, \ d^2=1 \ \rangle. \]

The P-positions are $\{a,b^2,bc,c^2,ac^3,ad,b^3d,cd,bc^2d \}$.  The pretending function appears in Figure \ref{071}.

\begin{figure}
\begin{center}
\[
\begin{array}{c|cccccc}
   & 1 & 2 & 3 & 4 & 5 & 6 \\ \hline
0+ & a & b & a & 1 & c & 1 \\
6+ & a & d & a & 1 & c & 1 \\
12+ & a & d & a & 1 & c & 1 \\
18+ & \ldots 
\end{array}\]
\end{center}
\caption{\label{071} The pretending function for misere play of {\bf 0.71}.}
\end{figure}

\subsection{Other quotients}

Hundreds more such solutions have been found amongst the octal games.  The forthcoming paper \cite{ps} includes
a census of such results.

\section{Computing presentations \& MisereSolver}

How are such solutions computed?  Aaron Siegel's recently developed Java program {\em MisereSolver} \cite{miseresolver}
will do it for you!  Some details on the algorithms used in {\em MisereSolver} are included in \cite{ps}.  Here,
we simply give a flavor of the some ideas underpinning it and how the software is used.

\subsection{Misere periodicity}

At the center of the Sprague-Grundy theory is the equation $G+G = 0$, 
which always holds for an arbitrary normal play combinatorial game $G$.
One consequence of $G+G=0$ is the equation
\[G+G+G = G,\]
in which all we've done is add $G$ to both sides.
In general, in normal play,
\begin{equation*}(k+2) \cdot G = k \cdot G. \end{equation*}
holds for every $k \geq 0$.

In {\em misere play}, the relation
\[(G + G) \ \rho \ \ \mathbbm{o} \]
happens to be true for $G= \mathbbm{o}$ and $G=\mathbbm{1}$, but beyond that, it is only seldom true for
occasional rule sets $\Gamma$ and positions $G$.
On the other hand, 
\[(G+G+G) \ \rho \ G\]
is {\em very often} true in misere play, and it is {\em always} true, for all $G$, if $\Gamma$ is a tame game.  And in {\em wild}
games $\Gamma$ for which the latter equation fails, often a weaker equation such as 
\begin{equation*}(G+G+G+G) \ \rho \  (G+G),\end{equation*}
is still valid, regardless of $G$.

These considerations suggest that a useful place to look for misere quotients is inside commutative monoids
having some (unknown) number of generators $x$ each satisfying a relation of the form
\[x^{k+2} = x^k\]
for each generator $x$ and some value of $k \geq 0$.

\subsection{Partial quotients for heap games}

A {\em heap game} is an impartial game $\Gamma$ whose rules can be expressed in terms of play on separated, non-interacting heaps of beans.
In constructing misere quotients for heap games, it's useful to introduce the {\em $n$th partial quotient}, which is just the indistinguishability
quotient of $\Gamma$ when all heaps are required to have $n$ or fewer beans.

\subsection{MisereSolver output of partial quotients}

Here is an (abbreviated) log of {\em MisereSolver} output of partial quotients for {\bf 0.123}, an octal game that is studied in great detail in \cite{ttw}.
In this output, monomial exponents have been juxtaposed with the generator names (so that $b^2c$, for example, appears as {\tt b2c}).
The program stops when it discovers the entire quotient---the partial quotients stabilize in a monoid of order 20, whose single-heap pretending function $\Phi$ is periodic of length 5.

\begin{verbatim}
C:\work>java -jar misere.jar 0.123
=== Normal Play Analysis of 0.123 ===
Max   : G(3) = 2
Period: 5 (5)
=== Misere Play Analysis of 0.123 ===
-- Presentation for 0.123 changed at heap 1 --
Size 2: TAME
P = {a}
Phi = 1 a 1
-- Presentation for 0.123 changed at heap 3 --
Size 6: TAME
P = {a,b2}
Phi = 1 a 1 b b a b2 1
-- Presentation for 0.123 changed at heap 8 --
Size 12: {a,b,c | a2=1,b4=b2,b2c=b3,c2=1}
P = {a,b2,ac}
Phi = 1 a 1 b b a b2 1 c
-- Presentation for 0.123 changed at heap 9 --
Size 20: {a,b,c,d | a2=1,b4=b2,b2c=b3,c2=1,b2d=d,cd=bd,d3=ad2}
P = {a,b2,ac,bd,d2}
Phi = 1 a 1 b b a d2 1 c d a d2 1 c d a d2 1 c d a d2 1
=== Misere Play Analysis Complete for 0.123 ===
Size 20: {a,b,c,d | a2=1,b4=b2,b2c=b3,c2=1,b2d=d,cd=bd,d3=ad2}
P = {a,b2,ac,bd,d2}
Phi = 1 a 1 b b a d2 1 c d a d2 1 c d a d2 1 c d a d2 1
Standard Form : 0.123
Normal Period : 5
Normal Ppd    : 5
Normal Max G  : G(3) = 2
Misere Period : 5
Misere Ppd    : 5
Quotient Order: 20
Heaps Computed: 22
Last Tame Heap: 7
\end{verbatim}

\subsection{Partial quotients and pretending functions}

Let's look more closely at the {\em MisereSolver} partial quotient output in order to 
illustrate some of the subtlety of misere quotient presentation calculation.

In Figure \ref{446}, we've shown three pretending functions for {\bf 0.123}.  The first is just
the normal play pretending function (ie, the nim-sequence) of the game, to heap six.  The second table shows
the corresponding misere pretending function for the partial quotient to heap size $6$, and the final table shows the initial portion of the pretending function 
for the entire game (taken over arbitrarily large heaps).

With these three tables in mind, consider the following question:

\[\mbox{When is } 4+4 \mbox{ indistinguishable from } 6 \mbox { in {\bf 0.123}}? \]

\begin{figure}[h]
Normal ${\bf 0.123}$

\begin{tabular}{c|ccccccccccc}
$n$ & 1 & 2 & 3 & 4 & 5 & 6 & 7 & 8 & 9 & 10 \\ \hline
$G(n)$ & $\ast 1 $ & $\ast 0 $ & $\ast 2$  & {$\ast 2$} & $\ast 1 $ & {$\ast 0$} & $\cdots $ & $\cdots$ & $\cdots$ & $\cdots$ \\ 
\end{tabular}

\vspace{0.2in}

Misere ${\bf 0.123}$ to heap 6:  $\langle a,b \ | \ a^2=1, \ b^3=b \rangle$, order 6

\begin{tabular}{c|cccccccc}
$n$    & 1    & 2    & 3    & 4    & 5  & 6  \\ \hline 
$\Phi(n)$  & $a$  & $1$  & $b$  & {$b$} & $a$  & {$b^2$}     \\
\end{tabular}

\vspace{0.2in}

Complete misere ${\bf 0.123}$ quotient, order 20
\[\langle a,b,c,d \ | \ a^2=1,\ b^4=b^2,\ b^2c=b^3,\ c^2=1,\ b^2d=d,\ cd=bd,\ d^3=ad^2 \rangle \]

\begin{tabular}{c|ccccccccccc}
$n$  & 1 & 2 & 3 & 4 & 5 & 6 & 7 & 8 & 9 & 10 \\ \hline 
$\Phi(n)$ & $a$  & $1$  & $b$  & {$b$}  & $a$  & {$d^2$}  & $1$  & $c$  & $d$  &  $\cdots$ \\  
\end{tabular}
\caption{\label{446} Iterative calculation of misere partial quotients differs in a fundamental way from
normal play nim-sequence calculation because sums at larger heap sizes (for example, 8+9) may distinguish between positions that previously
were indistinguishable at earlier partial quotients (eg, 4+4 and 6, to heap size six).}
\end{figure}

Let's answer the question.
In normal play (the top table), 4+4 is indistinguishable from 6 because 
\[G(4+4) = G(4)+G(4) = \ast 2 + \ast 2 = \ast 0 = G(6).\]  And in
the middle table, 4+4 is also indistinguishable from 6, since both sums evaluate to $b^2$.
But in the final table, 
\[\Phi(4+4) = \Phi(4)+\Phi(4) = b \cdot b = b^2 \neq d^2 = \Phi(6),\]  
ie, 4+4 can be distinguished from 6 in play of {\bf 0.123} {\em when no restriction is placed on the heap sizes}.
In fact, one verifies that the sum $8+9$, a position of type $cd$, distinguishes between 4+4 and 6 in {\bf 0.123}.

The fact that the values of partial misere pretending functions may change in this way, as larger heap sizes
are encountered, makes it highly desirable to carry out the calculations via computer programs that know how
to account for it.  

\subsection{Quotients from canonical forms}

In addition to computing quotients directly from the Guy-Smith code of octal games \cite{gs},
{\em MisereSolver} also can take as input the a canonical form of a misere game $G$.  It then computes
the indistinguishability quotient of its closure cl($G$).  This permits more general games than simply heap games to be analyzed.

\subsubsection{A coin-sliding game}
\label{coinslide}
For example, suppose we take $G=\{\mathbbm{2}_+, \mathbbm{o}\}$, a game listed in Figure \ref{conway-list}. 
In the output script below, {\em MisereSolver} calculates that the indistinguishability quotient of cl($G$) 
is a monoid of order 14 with four P-position types:

\begin{figure}
\includegraphics{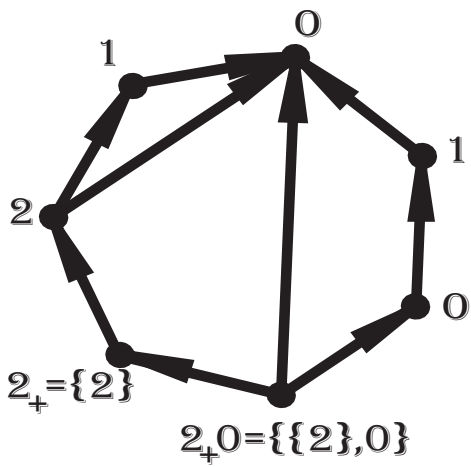}
\caption{\label{hept} Misere coin-sliding on a directed heptagon with two additional edges.   An arbitrary number of coins are placed at the vertices, and two players take turns
sliding a single coin along a single directed edge.  Play ends when the final coin reaches the topmost (sink) node (labelled $\mathbbm{o}$).  Whoever
makes the last move loses the game.  The associated indistinguishability quotient is a commutative monoid of order 14 with presentation
\[\langle \ a, \ b, \ c \ | \ a^2=1,\ b^3=b,\ b^2c =c,\ c^3=ac^2 \ \rangle \mbox{\ \ \ \ \ \ \ \ \ \ \ \ \ \ \ \ }\] and P-positions $\{\ a,\ b^2,\ bc,\ c^2\}$.  See section \ref{coinslide} and Figure \ref{h-corr}.}.   
\end{figure}

{\small
\begin{verbatim}
-- Presentation for 2+0 changed at heap 1 --
Size 2: TAME
P = {a}
Phi = 1 a
-- Presentation for 2+0 changed at heap 2 --
Size 6: TAME
P = {a,b2}
Phi = 1 a b b2
-- Presentation for 2+0 changed at heap 4 --
Size 14: {a,b,c | a2=1,b3=b,b2c=c,c3=ac2}
P = {a,b2,bc,c2}
Phi = 1 a b c2 c
\end{verbatim}
}

Figure \ref{hept} shows a coin-sliding game that can be played perfectly using this information.  Figure \ref{h-corr} shows how the
canonical forms at each vertex correspond to elements of the misere quotient.

\begin{figure}[h]
\[
\begin{array}{c|cccccc}
\mbox{Canonical form} & \mathbbm{o} & \mathbbm{1} & \mathbbm{2} & \mathbbm{2}_+ & \{\mathbbm{2}_+,\mathbbm{o}\}    \\ \\ \hline
\\
\mbox{Quotient element} & 1 & a & b & c^2 & c
\\
\end{array}\]
\caption{\label{h-corr} Assignment of single-coin positions in the heptagon game to misere quotients elements.}
\end{figure}

\section{Outlook}

At the time of this writing (December 2005), the indistinguishability quotient construction is only one year old.
Several aspects of the theory are ripe for further development, and the misere versions of 
many impartial games with complete normal play solutions remain to be investigated.  We have space only to describe a few of the
many interesting topics for further investigation.

\subsection{Infinite quotients}  Misere quotients are not always finite.  Today, it frequently happens that {\em MisereSolver}
will ``hang'' at a particular heap size as it discovers more and more distinguishable position types.  Is it possible to
improve upon this behavior and discover algorithms that can handle infinite misere quotients?  

\subsubsection{Dawson's chess}  
\label{dawson} One important game that seems to have an infinite misere quotient is Dawson's Chess.
In the equivalent form {\bf 0.07}, (called Dawson's Kayles), Aaron Siegel \cite{ps} found that the order of its misere partial quotients ${\mathcal Q}$ grows as indicated
in Figure \ref{dawson-qs}:

\begin{figure}[h]
\[
\begin{array}{c|cccccccc}
\mbox{Heap size} & 24 & 26 & 29  & 30 & 31 & 32 & 33 & 34   \\ \\ \hline
\\
\mbox{$|Q|$} & 24 & 144 & 176 & 360 & 520 & 552 & 638 & \infty (?)
\\
\end{array}\]
\caption{\label{dawson-qs} Is {\bf 0.07} infinite at heap 34?}
\end{figure}

Since Redei's Theorem (see \cite{ttw} for discussion and additional references) asserts that a finitely generated commutative monoid is
always finitely presentable, the object being sought in Figure \ref{dawson-qs} (the misere quotient presentation to heap size 34) certainly
exists, although it most likely has a complicated structure of P- and N-positions.  New ideas are needed here.

\subsubsection{Infinite, but not at bounded heap sizes}

Other games seemingly exhibit infinite behavior, but appear to have finite order (rather than simply finitely presentable) partial quotients at all heap sizes.
One example is {\bf .54}, which shows considerable structure in the partial misere quotients output by {\em MisereSolver}.  Progress
on this game would resolve difficulties with
an incorrect solution of this game that appears in the otherwise excellent paper \cite{a3}.  Siegel
calls this behavior {\em algebraic periodicity}.

\subsection{Classification problem}

The {\em misere quotient classification problem} asks for an enumeration of the possible nonisomorphic misere quotients at each order
$2k$, and a better understanding of the category of commutative monoids that arise as misere quotients\footnote{It can be shown that a finite misere quotient has even order \cite{ps}.}. 
Preliminary computations by Aaron Siegel suggest that the number of nonisomorphic misere quotients grows as follows:

\begin{figure}[h]
\[
\begin{array}{c|ccccccccc}
\mbox{Order} & 2 & 4 & 6  & 8 & 10 & 12 
\\ \hline
\mbox{\# quotients} & 1 & 0 & 1 & 1 & 1? & 6? 
\\
\end{array}\]
\caption{\label{dawson-qs} Conjectured number of nonisomorphic misere quotients at small orders.}
\end{figure}

Evidently misere quotients are far from general commutative 
semigroups---by comparison, the number of nonisomorphic 
commutative semigroups at orders 4, 6, and 8  are already 58, 2143, and 221805, respectively (\cite{gril}, pg 2).

\subsection{Relation between normal and misere play quotients}

If a misere quotient is finite, does each of its elements $x$ necessarily satisfy a relation of the form $x^{k+2} = x^{k}$, for
some $k \geq 0$?  The question is closely related to the structure of maximal subgroups inside misere finite quotients.
Is every maximal subgroup of the form $(\mathbb{Z}_2)^m$, for some $m$?

At the June 2005 Banff conference on combinatorial games, the author conjectured that an octal game, if
misere periodic, had a periodic normal play nim sequence with the two periods (normal and misere) equal.  Then Aaron Siegel
pointed out that {\bf 0.241}, with normal period two, has misere period 10.  Must the normal period length {\em divide} the misere one,
if both are periodic?

\subsection{Quaternary bounties}

Again at the Banff conference, the author distributed the list of wild misere quaternary games in Figure \ref{quatwild4}.

\begin{figure}[h]
\begin{center}
\begin{tabular}{llll}
$({\bf .0122}, 1^{20}, 12)$ & $({\bf .0123}, 1^{20}, 12)$ & $({\bf .1023}, 2^{1420}, 11)$ & $({\bf .1032}, 2^{1420},12)$ \\
$({\bf .1033}, 1^{20}, 11)$ & $({\bf .1231}, 2^{1420}, 8)$ & $({\bf .1232}, 2^{1420}, 9)$ & $({\bf .1233}, 2^{1420}, 9)$ \\
$({\bf .1321}, 2^{1420}, 9)$ & $({\bf .1323},2^{1420},10)$ & $({\bf .1331}, 1^{20}, 8)$ & $({\bf .2012}, 1^{20}, 5)$ \\
$({\bf .2112}, 1^{20},5)$ & $({\bf .3101},1^{20},4)$ & $({\bf .3102},0^{20},5)$ & $({\bf .3103}, 1^{20}, 4)$ \\ 
$({\bf .3112}, 2^{1420}, 7)$ & $({\bf .3122}, 2^{1420},4)$ & $({\bf .3123}, 1^{31}, 6)$ & $({\bf .3131}, 2^{1420}, 6)$ \\
$({\bf .3312}, 2^{1420}, 5)$ & & &
\end{tabular}
\caption{\label{quatwild4} The twenty-one wild four-digit quaternary games (with first wild genus value \& corresponding heap size)}
\end{center}
\end{figure}

The author offered a bounty of \$25 dollars/game to the first person to exhibit the misere indistinguishability quotient and pretending function of the games in the list.
Aaron Siegel swept up 17 of the bounties \cite{ps}, but {\bf .3102},  {\bf .3122}, {\bf .3123}, and {\bf .3312} are still open.

\subsection{Misere sprouts endgames}  Misere Sprouts (see \cite{ww}, 2nd edition, Vol III) is perhaps the only 
misere combinatorial game that is played competitively in an organized forum, the {\em World Game of Sprouts Association}.
It would be interesting to assemble a database of misere sprout endgames and compute the indistinguishability quotient of their misere addition.

\subsection{The misere mex mystery}  In normal play game computations for heap games, the {\em mex rule} allows the computation of the heap $n+1$
nim-heap equivalent from the equivalents at heaps of size $n$ and smaller.  The {\em misere mex mystery} asks for the analogue of the normal play mex rule, in misere play.  It is evidently closely related to the partial quotient computations performed by {\em MisereSolver}.

\subsection{Commutative algebra}

A beginning at application of theoretical results on commutative monoids to misere quotients was begun in \cite{ttw}.  What more can be said?

\section{Appendix: The genus theory}
\label{gthm}
We summarize Conway's {\em genus theory}, first described in chapter 12 of \cite{onag}, and used extensively in {\em Winning Ways}.
It describes a method for calculating whether all the positions of particular game $\Gamma$ are tame, and how to give a complete
analysis of $\Gamma$, if so.  The genus
theory assigns to each position $G$ a particular symbol
\begin{equation}
\label{genuseq}
\mbox{genus}(G)= {\sf G}^{\ast}(G)= g^{g_0g_1g_2\cdots}.
\end{equation}
where the $g$ and the $g_i$'s are always nonnegative integers.
We'll define this genus value precisely and illustrate how to calculate genus values for some example games $G$, below.

To look at this in more detail, we need some preliminary definitions before giving definition of genus values.

\subsection{Grundy numbers}

Let $\ast k$ represent the nim heap of size $k$. The {\em Grundy number} 
(or {\em nim value}) of an impartial game position $G$
is the unique number $k$ such that $G + \ast k$  is a second-player win.
Because Grundy numbers may be defined relative to normal or misere play,
we distinguish between the {\em normal play Grundy number} ${\sf G}^{+}(G)$
and its counterpart ${\sf G}^{-}(G)$, the {\em misere Grundy number}.
  
In normal play,  Grundy numbers can be calculated using the rules
 ${\sf G}^{+}(0) = 0$, and otherwise, ${\sf G}^{+}(G)$ is the 
least number (from 0,1,2, \ldots) that is {\em not} the Grundy
number of an option of $G$ (the so-called {\em minimal excludant}, 
or {\em mex}).

When normal play is in effect, every game with Grundy number
${\sf G}^{+}(G)= k$ can be thought of as the nim
heap $\ast k$.  No information about best play of the game 
is lost by assuming that $G$ is in fact
precisely the nim heap of size $k$.   Moreover, in normal play, the Grundy number of
a sum is just the nim-sum of the Grundy numbers of the summands.

The misere Grundy number 
is also simple to define (see \cite{onag}, pg 140, bottom):

{\sf
\begin{quotation}
 ${\sf G}^{-}(0) = 1$.  Otherwise, ${\sf G}^{-}(G)$ is the least number (from 0,1,2, \ldots)
which is not the ${\sf G}^{-}$-value of any option of $G$.  Notice that this is just like
the ordinary ``mex'' rule for computing ${\sf G}^{+}$, except that we have 
 ${\sf G}^{-}(0) = 1,$ and ${\sf G}^{+}(0) = 0$.  
\end{quotation}
}

Misere P-positions are precisely those whose first genus exponent is 0.

\subsection{Indistinguishability vs misere Grundy numbers}

When misere play is in effect, Grundy numbers can still be defined---as we've already 
said---but many {\em distinguishable} games are assigned the {\em same} Grundy number, and the outcome
of a sum is {\em not} determined by Grundy numbers of the summands.  These unfortunate
facts lead directly to the apparent great complexity of many misere analyses.  

Here is the definition of the genus, directly from \cite{onag}, now at the bottom of page 141:

{\sf
\begin{quotation}
In the analysis of many games, we need even more information than is provided by either
of these values [${\sf G}^{+}$ and ${\sf G}^{-}$],  
and so we shall define a more complicated symbol that we call the 
${\sf G}^{*}$-value, [or {\em genus}\ ], ${\sf G}^{*}(G)$.  This is the symbol

\[g^{g_0g_1g_2\cdots}\]
where 
\begin{eqnarray*}
g & = & {\sf G}^{+}(G) \\
g_0 & = & {\sf G}^{-}(G) \\
g_1 & = & {\sf G}^{-}(G + \mathbbm{2}) \\
g_2 & = & {\sf G}^{-}(G+ \mathbbm{2}+ \mathbbm{2}) \\
\ldots & = & \ldots \\
\end{eqnarray*}
where in general $g_n$ is the ${\sf G}^{-}$-value of the sum of $G$ with $n$ other games all
equal to [the nim-heap of size] 2.
\end{quotation}
}  

At first sight, the genus symbol looks to be an potentially infinitely long symbol in its 
``exponent.''  In practice, it can be shown that the $g_i$'s always fall into
an eventual period two pattern.  By convention, a genus symbol is 
written down with a finite exponent with the understanding that its final two
values repeat indefinitely.  

The only genus values that arise in misere Nim are the {\em tame genera}

\begin{figure}[h]
\[\begin{array}{c}
\underbrace{0^{120}, 1^{031}}_{\stackrel{\mbox{Genera of normal play $\ast 0$ (resp, $\ast 1$) Nim}}{\mbox{positions involving nim heaps of size 1 only;}}} \ 
\\ 
\\ and
\\
\\
\underbrace{{0^{02},1^{13}, 2^{20}, 3^{31}, 4^{46}, \cdots, n^{n(n \oplus 2)}}, \cdots}_{\stackrel{\mbox{Genera of $\ast n$ normal-play Nim positions}}{\mbox{involving at least one nim heap of size } \geq 2.}}

\end{array}\]
\caption{\label{genusNim} Correspondence between normal play nim positions and tame genera.}
\end{figure}
The value of the genus theory lies in the following theorem (cf \cite{onag}, Theorem 73):

{\sf
\begin{quotation}
{\bf Theorem}: If all the positions of some game $\Gamma$ have tame genera, the genus of a sum $G+H$ can be computed by replacing
the summands by Nim-positions of the same genus values, and taking the genus value of the resulting sum.
\end{quotation}
}

In order to apply the theorem to analyze a tame game $\Gamma$, a person needs to know several things:

\begin{enumerate}
\item How to compute genus symbols for positions $G$ of a game $\Gamma$;
\item That every position of the game $\Gamma$ does have a tame genus;
\item The correspondence between the tame genera and Nim positions.
\end{enumerate}
We've already given the correspondence between normal-play Nim positions and their misere genus values, in Figure (\ref{genusNim}).
We'll defer the most complicated part---how to compute genera, and verify that they're all tame---to the next section.

The addition rule for tame genera is not complicated. The first two 
symbols have the $\mathbb{Z}_2$ addition
\begin{eqnarray*}
 0^{120} + 0^{120} & = & 0^{120} \\
 0^{120} + 1^{031} & = & 1^{031} \\
 1^{031} + 1^{031} & = & 0^{120} 
\end{eqnarray*}

Two positions with genus symbols of the form $n^{n(n \oplus 2)}$ add just like Nim heaps of $\ast n$. For example,
\[2^{20} + 3^{31} = 1^{13}. \]
The symbol $0^{120}$ adds like an identity, for example:
\[4^{46} + 0^{120} = 4^{46}. \]
When $1^{031}$ is added to a $n^{n(n \oplus 2)}$, it acts like $1^{13}$: 
\[4^{46} + 1^{031} = 5^{57}. \]
It has to emphasized that these rules work {\em only if all positions in play of $\Gamma$ are known to have tame genus values}.
If, on the other hand, even a single position in a game $\Gamma$ does {\em not} have a tame
genus, the game is wild and {\em nothing can be said in general about the addition of tame genera}.

\subsection{Genus calculation in octal game {\bf 0.123}}

Let's press on with the genus theory, illustrating it in an example game, and 
keeping in mind the end of Chapter 13 in \cite{ww}:

{\sf
\begin{quotation}
{\em The misere theory of impartial games is the last and most complicated
theory in this book.  Congratulations if you've followed us so far...}
\end{quotation}
}

Genus computations, and the nature of the conclusions that can be
drawn from them, are what makes Chapter 13 in 
{\em Winning Ways} complicated.  In this
section we illustrate genus computations by using them to 
initiate the analysis of a particular wild octal game ({\bf 0.123}).
Because the game {\bf 0.123} is wild, the genus theory will {\em not}
lead to a complete analysis of it.  A complete analysis can nevertheless be obtained
via the indistinguishability quotient construction; for details, see \cite{ttw}.

The octal game {\bf 0.123} can be played with counters arranged in heaps.
Two players take turns removing one, two or three counters from a heap,
subject to the following additional conditions:
\begin{enumerate}
\item Three counters may be removed from any heap;
\item Two counters may be removed from a heap, but only if it has more than two counters; and
\item One counter may be removed only if it is the only counter in that heap.
\end{enumerate}

\subsubsection{Normal play of {\bf 0.123}} 
The nim sequence of {\bf 0.123}\footnote{See {\it Winning Ways}, Chapter 4, pg 97, ``Other Take-Away Games;'' also Table 7(b),
pg 104.} is periodic of length 5, beginning at heap 5.  See Figure \ref{123normal}.

\begin{center}
\begin{figure}
\begin{tabular}{c|ccccc}
+   & 1 & 2 & 3 & 4 & 5 \\ \hline
0+ & 1 & 0 & 2 & 2 & 1 \\
5+ & 0 & 0 & 2 & 1 & 1 \\
10+ & 0 & 0 & 2 & 1 & 1 \\
15+ & $\cdots$ & & & &
\end{tabular}
\caption{\label{123normal}  Normal play nim values of {\bf 0.123}}
\end{figure}
\end{center}

\subsubsection{Misere play genus computations for {\bf 0.123}} 

We exhibit single-heap genus values of {\bf 0.123} in Figure \ref{123genera}.  It's possible to prove that
this sequence is also periodic of length 5.  However, a periodic genus sequence is not the same thing as a complete
misere analysis. Let's see what happens instead.

\begin{figure}[h]
\begin{tabular}{c|ccccc}
+   & 1 & 2 & 3 & 4 & 5 \\ \hline
0+ & $1^{031}$ & $0^{120}$ & $2^{20}$ & $2^{20}$ & $1^{031}$ \\
5+ & $0^{02}$ & $0^{120}$ & $2^{1420}$ & $1^{20}$ & $1^{031}$ \\
10+ & $0^{02}$ & $0^{120}$ & $2^{1420}$ & $1^{20}$ & $1^{031}$ \\
15+ & $\cdots$ & & & &
\end{tabular}
\\
\caption{\label{123genera}{\sf G}*-values of {\bf 0.123}}
\end{figure}

There are some tame genus symbols in Figure \ref{123genera}.  They are
\begin{eqnarray*}
0 & = & 0^{1202020\cdots} = 0^{120} \\
1 & = & 1^{0313131\cdots} = 1^{031} \\
2 & = & 2^{2020202\cdots} = 2^{20}
\end{eqnarray*}

Despite the presence of these tame genera, the game is still wild---the first wild genus value, $2^{1420}$,
occurs at heap 8.  Conway's Theorem 73 on tame games therefore does {\em not} apply, since it requires {\em all} positions
to have tame genera in order for the game to be treated as misere Nim.  We can say nothing about how
genera add---even the tame genera---without examining the game more closely.

Here's what we can (and cannot) do with Figure \ref{123genera}.

\subsubsection{Single heaps}

We {\em can} determine the outcome class of {\em single-heap} {\bf 0.123} positions.
The first superscript in a heap's 
genus symbol is 0 if and only if that heap size
is a $P$-position.  The single heap $P$-positions of {\bf 0.123} therefore occur at heap sizes 
\[ 1, 5, 6, 10, 11, 15, 16, 20, 21, \ldots \]  
For example, the genus of the heap of size 7 has its first superscript = 1.
It is therefore an $N$-position.  The winning move is $7 \rightarrow 5$.

\subsubsection{Multiple heaps}
 We {\em cannot} immediately determine the outcome class of {\em multiple-heap}
{\bf 0.123} positions using Figure \ref{123genera}.  However, Figure \ref{123genera} does provide a basis for
investigating multiheap positions.  
For example, Figure \ref{123multiple} is a table that shows the genera of two-heap
positions up to heap size nine.

\begin{center}
\begin{figure}
\begin{tabular}{c|lllllllll}
+ & $h_{1}$ & $h_2$ & $h_3$ & $h_4$ & $h_5$ & $h_6$ & $h_7$ & $h_8$ & $h_9$
\\ \hline
$h_1$ & $0^{120}$ & $1^{031}$ & $3^{31}$ & $3^{31}$ & $0^{120}$ & $1^{13}$ & $1^{031}$ & $3^{0531}$ & $0^{31}$ \\
$h_2$ &     & $0^{120}$ & $2^{20}$ & $2^{20}$ & $1^{031}$ & $0^{02}$  & $0^{120}$ & $2^{1420}$ & $1^{20}$ \\
$h_3$ &     &     & $0^{02}$ & $0^{02}$ & $3^{31}$ & $2^{20}$ & $2^{20}$ & $0^{420}$ & $3^{02}$ \\
$h_4$ &     &     &     & $0^{02}$ & $3^{31}$ & $2^{20}$ & $2^{20}$  & $0^{420}$ & $3^{02}$ \\
$h_5$ &     &     &     &       & $0^{120}$ & $1^{13}$ & $1^{031}$ & $3^{0531}$ & $0^{31}$ \\
$h_6$ &     &     &     &       &     & $0^{02}$ & $0^{02}$ & $2^{20}$ & $1^{13}$ \\
$h_7$ &     &     &     &       &     &       & $0^{120}$ & $2^{1420}$ & $1^{20}$ \\
$h_8$ &     &     &     &       &     &       &     & $0^{120}$ & $3^{02}$ \\
$h_9$ &     &     &     &       &     &       &     &          & $0^{02}$ \\

\end{tabular}
\\
\caption{\label{123multiple} Some genus values of games $h_i + h_j$ in {\bf 0.123}.}
\end{figure}
\end{center}

\subsection{Genus calculation algorithm} Here's how the genus of a particular sum $G=h_8 + h_5$ was computed from
the earlier single-heap values in Figure \ref{123genera}.
First, we rewrote genus(G) in terms of its options:
\[{\rm genus}(G) = {\rm genus}(h_8 + h_5) 
= {\rm genus}(\{h_6 + h_5, h_5 + h_5, h_8+h_3, h_8+h_2\}) \]
The genus of a non-empty game $G = \{A, B, \cdots\}$ can be 
calculated
from the genus of its options
$A, B, \ldots$
using the 
{\em mex-with-carrying algorithm} 
($\diamond$ symbols represent positions with no carry):

\begin{eqnarray*}
        {\rm carry } (\gamma)              & = & {\diamond}^{\diamond05313} \\
        {\rm carry } (\gamma \oplus 1)         & = & {\diamond}^{\diamond14202} \\
{\rm genus}(h_6 + h_5) & = & 1^{131313\ldots} \\
{\rm genus}(h_5 + h_5) & = & 0^{120202\ldots} \\
{\rm genus}(h_8 + h_3) & = & 0^{420202\ldots} \\
{\rm genus}(h_8 + h_2) & = & \underline{2^{142020\ldots}} \\
{\rm genus}(G)       & = & 3^{053131\ldots} 
\end{eqnarray*}
\vspace{0.1in}

The result ${\rm genus}(G) =  3^{053131\ldots} =  3^{0531}$
was computed columnwise, working from left to right.  First, the ``base'' 
and ``first superscript'' results
\[{\sf G}^{+}(G) = {\rm mex}(\{1,0,0,2\}) = 3\]
and
\[{\sf G}^{-}(G) = {\rm mex}(\{1,1,4,1\}) = 0\]
were computed from the corresponding four positions in each option of $G$,
with no carries present.  The ``carry out'' is then $\gamma = 0$. 
The second superscript result
\[{\sf G}^{-}(G+ \ast 2) = {\rm mex}(\{3,2,2,4,{\bf 0},{\bf 1}\}) = 5\]
involved a similar computation, but with two {\em carry values}
\[\{\gamma,\gamma \oplus 1 \} = \{0,1\}.\]
thrown into the mex calculation (they're shown in bold).    
See the more complete description of this 
algorithm in the section titled {\em ``But What if They're Wild?''\ asks the Bad Child} (\cite{ww}, page 410).  
It's also illustrated on pg 143 in \cite{onag}.

\end{document}